\documentclass{article}

\usepackage{graphicx}
\usepackage{subfigure}
\usepackage{multirow}
\usepackage{listings}
\usepackage{amsmath}
\usepackage{amsfonts}
\usepackage{caption}
\usepackage[margin=1in]{geometry}
\usepackage{soul,color}
\usepackage{cite}

\begin{document}
\title{Randomized Sch\"{u}tzenberger's jeu de taquin and approximate calculation of co-transition probabilities of a central Markov process on the 3D Young graph
}

\author{Vasilii Duzhin, Nikolay Vassiliev}


\date{}
\maketitle

\begin{abstract}

There exists a well-known hook-length formula for calculating the dimensions of 2D Young diagrams. Unfortunately, the analogous formula for 3D case is unknown. We introduce an approach for calculating the estimations of dimensions of three-dimensional Young diagrams also known as plane partitions. 
The most difficult part of this task is the calculation of co-transition probabilities for a central Markov process. 
We propose an algorithm for approximate calculation of these probabilities. It generates numerous random paths to a given diagram. In case the generated paths are uniformly distributed, the proportion of paths passing through a certain branch gives us an approximate value of the co-transition probability. As our numerical experiments show, the random generator based on the
randomized variant of the  
 Sch\"{u}tzenberger transformation allows to obtain accurate values of co-transition probabilities.
Also a method to construct 3D Young diagrams with large dimensions is proposed.
\end{abstract}

\section{Introduction}

A 2D Young tableau of a size $n$ is a 2D Young diagram of $n$ boxes filled by integers 1..$n$ in such a way that numbers grow along the coordinate axes. Obviously, it can be generalized to 3D case: a 3D Young tableau is a 3D Young diagram filled by integers 1..$n$ which grow along $x$, $y$ and $z$ axes.
The Sch\"{u}tzenberger transformation on Young tableaux, also known as ``jeu de taquin'', was introduced in Sch\"{u}tzenberger's paper \cite{schu}. This transformation was applied to solve different problems of enumerative combinatorics and representation theory of symmetric groups. Particularly, it can be used to calculate the Littlewood--Richardson coefficients \cite{fomin}. 

Different Markov processes can be defined on the Young graph by assigning probabilities to it edges.
There exists the central Plancherel process on the two-dimensional Young graph. Because of its centrality, the probabilities of all paths from the root to a diagram of a certain shape are equal. So, the conditional Plancherel probabilities give the uniform distribution on the set of Young tableaux of a certain Young diagram. Therefore, we could use the conditional Plancherel probabilities for a uniform random generator of 2D Young tableaux of a fixed shape.

However, this approach is not applicable in the three-dimensional case because the central process on the 3D Young graph is unknown. By this reason, this method could not be used to generate random uniformly distributed 3D tableaux. Another possible approach to generate such tableaux is to efficiently enumerate tableaux of a certain shape. By efficient enumeration we mean a bijection of a set of all 3D tableaux of a fixed shape with a set of integers $[1..n] \in \mathbb{N}$ such that a tableau can be quickly restored by its number. But such an enumeration can be implemented only for tableaux of tens of thousands boxes even in 2D case. The complexity of enumeration of 3D Young tableaux is significantly higher, so such an enumeration cannot be practically used in numerical experiments involving large tableaux.

The goal of this work is to obtain the random generator of three-dimensional Young tableaux of the same shape which could be used for approximate calculation of co-transition probabilities of the central process. We propose the solution of this problem based on a special randomized modification of Sch\"{u}tzenberger's jeu de taquin. Also we describe an algorithm of constructing 3D Young diagrams with large dimensions based on this modification.

\section{Sch\"{u}tzenberger transformation}
\subsection{Standard Sch\"{u}tzenberger transformation}
We consider the Sch\"{u}tzenberger transformation on two- and three- dimensional Young tableaux. The Sch\"{u}tzenberger transformation converts a Young tableau of size $n$ to another Young tableau of size $n-1$.
At the beginning, the first box of a source tableau is being removed. Then, the box with a smaller number is being selected among top neighbouring and right neighbouring boxes. The selected box is then being shifted to the position of the removed box. A newly formed empty box is being filled by the neighbouring box using the same rule. This process continues until the front of the diagram is reached. Finally, all the numbers in boxes are decremented by one.

\begin{figure}[ht!]
        \centering
        \subfigure[]{
        \includegraphics[scale=0.06]{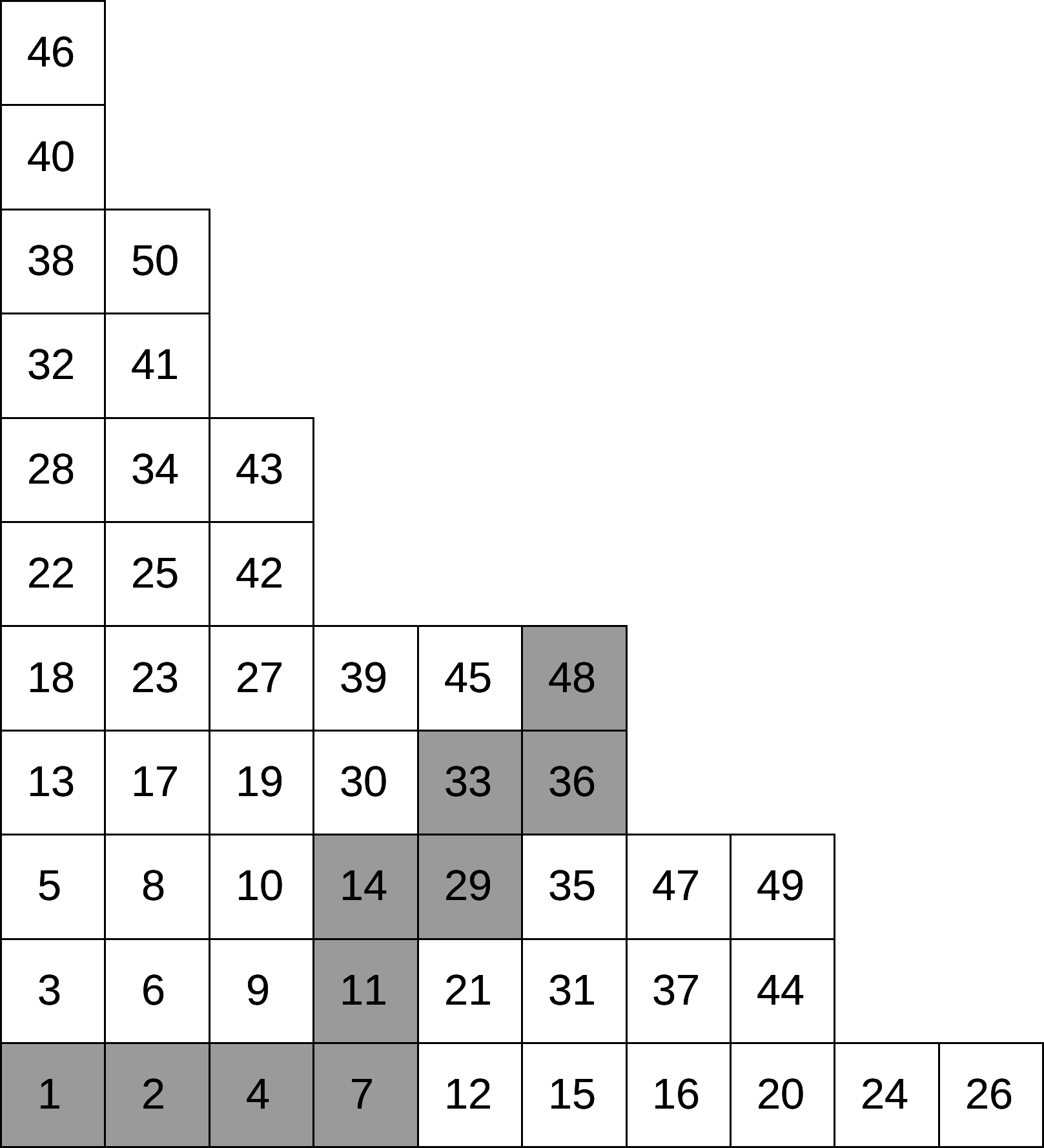}
        }
        \subfigure[]{
        \includegraphics[scale=0.06]{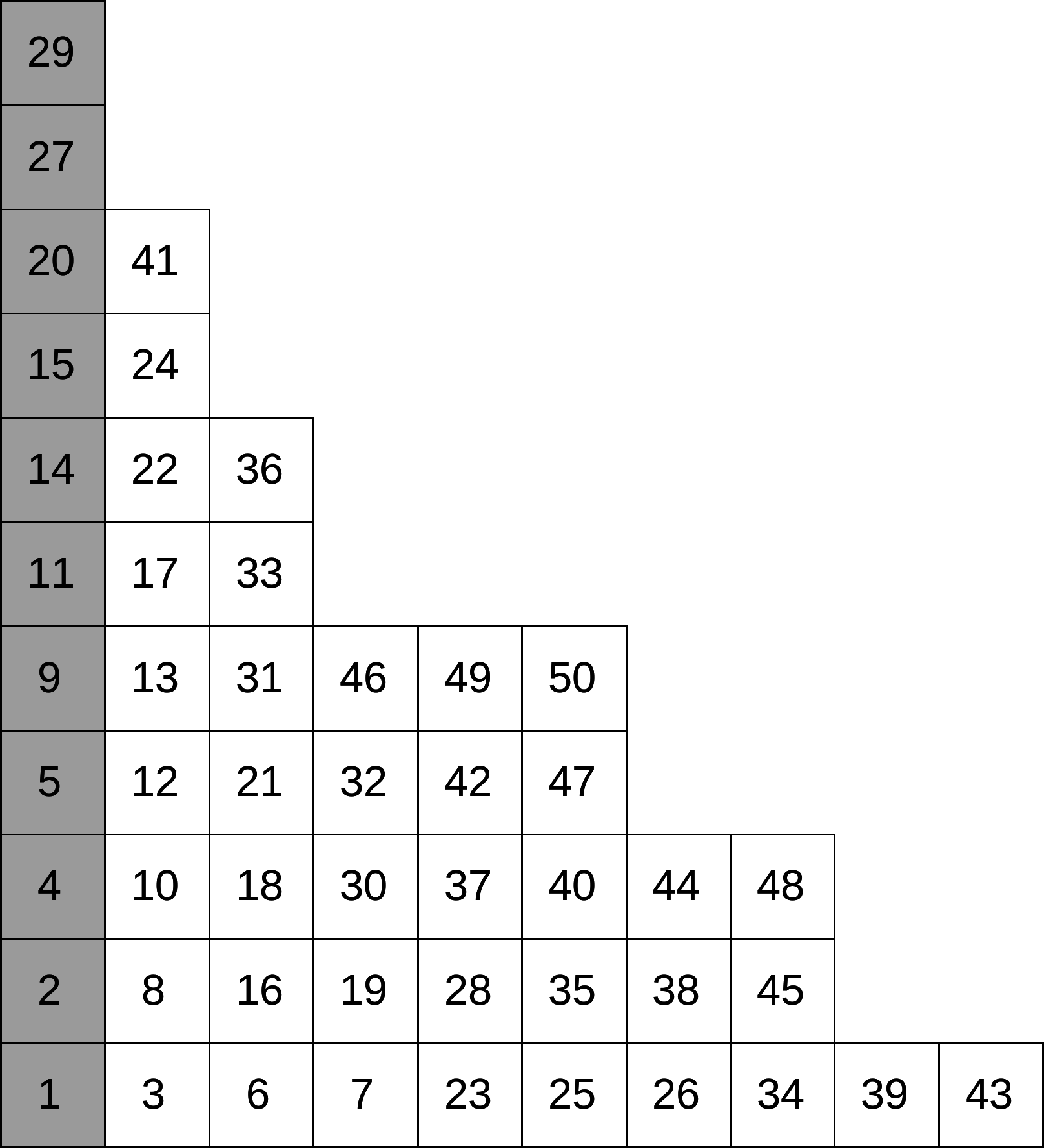}
        }
        \subfigure[]{
        \includegraphics[scale=0.06]{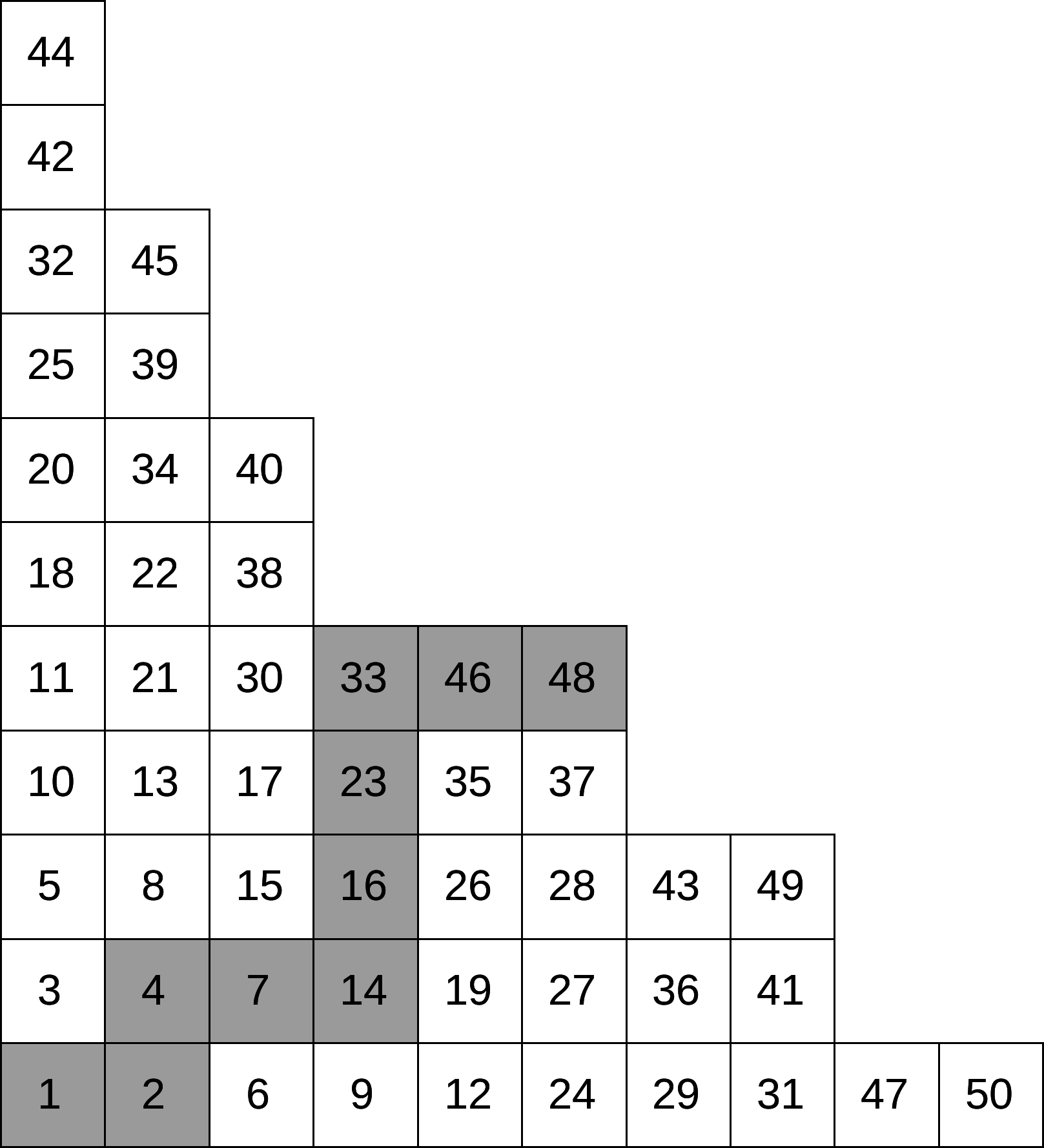}
        }
        \caption{Three different Young tableaux of the same shape and their nerves.}
        \label{fig:examples_schu}
\end{figure}

The sequence of the shifted boxes forms so-called \textit{jeu de taquin path} \cite{rom_snia} or \textit{Sch\"{u}tzenberger path}. A.~M.~Vershik suggests to call it \textit{a nerve}. The Sch\"{u}tzenberger paths are the paths in Pascal graphs: $\mathbb{Z}^2_{+}$ or $\mathbb{Z}^3_{+}$ in 2D and 3D cases, respectively.
Fig.~\ref{fig:examples_schu} shows some examples of Young tableaux and their nerves.

A.~M.~Vershik has noticed~\cite{vershik_arxiv} that jeu de taquin algorithm can be applied not only to the Young tableaux of an arbitrary dimension, but for a very broad class of partially ordered sets. In this case the Sch\"{u}tzenberger transformation works on ascendant sequences of decreasing ideals of a corresponding poset.
Particularly, the technique based on jeu de taquin can be used in a wide class of graded graphs.

In 2D case, many important properties of Sch\"{u}tzenberger transformation are connected closely to Plancherel measure. Particularly, it was proved in \cite{rom_snia} that a Sch\"{u}tzenberger path in two-dimensional case have a limit value of angle with a probability 1 relative to the Plancherel measure.

It would be natural to assume that 3D Sch\"{u}tzenberger transformation is also related to still unknown central measure in the 3D Young graph. Such a measure would be a 3D generalization of the Plancherel measure.

\subsection{Shape-preserving Sch\"{u}tzenberger transformation}
In this modification of the classic Sch\"{u}tzenberger transformation, we add an extra box in the position of the last shifted box. Thereby, this transformation does not change the shape of the diagram. Also the transformation becomes reversible, i.e. it establishes a bijection on the paths to a fixed vertex of the Young graph. This bijection splits into multiple cycles and computational experiments show that in most cases the number of cycles is relatively large. It means that the iterations of the shape-preserving Sch\"{u}tzenberger transformation do not give all possible Young tableaux of a certain shape. So, in order to build a uniform random generator of Young tableaux, we propose another modification of the transformation.

\subsection{Randomized Sch\"{u}tzenberger transformation}

Here we consider a randomized modification of the classic Sch\"{u}tzenberger transformation. Besides adding the last box of a jeu de taquin path to Young tableau in order to keep its original shape, we also choose the initial part of a tableau (or, equivalently, a path in the Young graph) randomly. 

\begin{figure}[ht]
        \centering
        \subfigure[]{
        \includegraphics[scale=0.12]{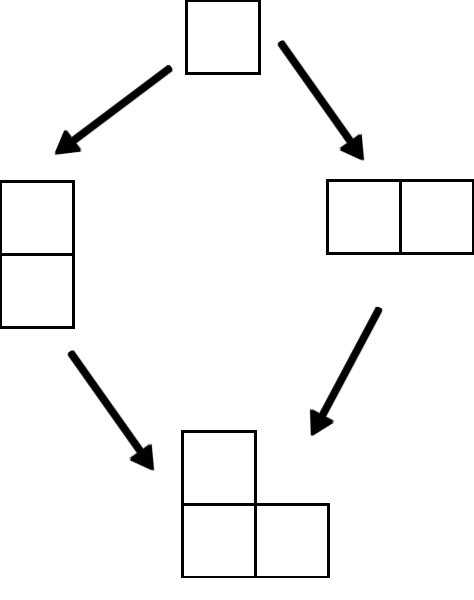}
        }
	\vrule
        \subfigure[]{
        \includegraphics[scale=0.2]{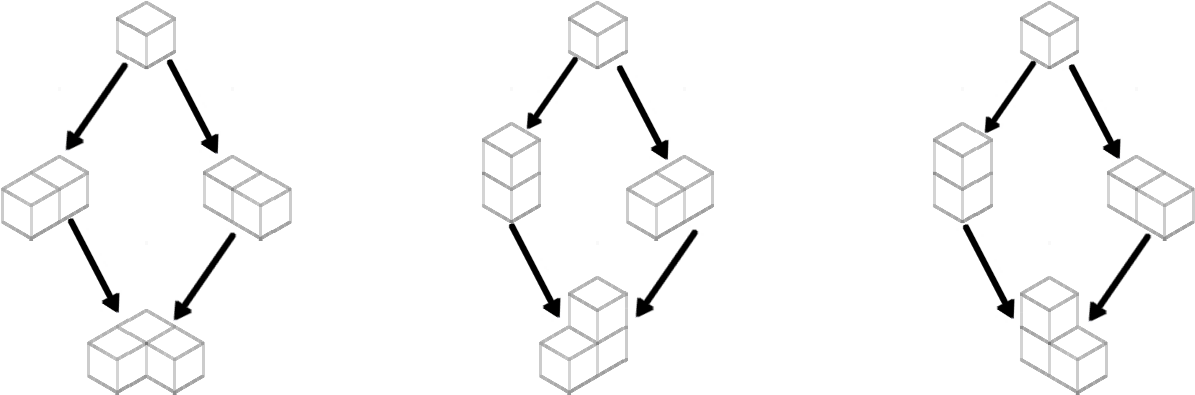}
        }
        \caption{Possible paths of randomization for (a): 2D Young tableaux, (b): 3D Young tableaux.}
        \label{fig:randomization}
\end{figure}

There exists the only two-dimensional diagram of size 3 with two incoming paths. In case the initial path passes through that diagram, with a probability of 1/2 we either keep the original path unchanged or switch to another one. Both paths are shown in Fig.~\ref{fig:randomization} (a).

There exists 3 different 3D Young diagrams of size 3 and dimension~2, see Fig.~\ref{fig:randomization} (b). If a tableau passes through one of these diagrams, the initial part of a path is either kept as it is or changed to another one with a probability of 1/2.

Note that the diagrams in Fig.~\ref{fig:randomization} are the minimum possible diagrams with dimension 2.
The results of numerical experiments discussed in Sec.~\ref{sec:numexp} allow to suggest that the iterations of the randomized Sch\"{u}tzenberger transformation generate almost uniform distribution on the paths to a diagram. 

\section{Fast Sch\"{u}tzenberger transformation}

We propose the following algorithm for the implementation of the Sch\"{u}tzenberger transformation on 2D and 3D Young tableaux. We present Young tableaux as one-dimensional arrays of coordinate pairs of added boxes. This allows to store tableaux in a compact way wherein it makes possible to restore the corresponding plane partition and to calculate the coordinates of boxes which can be added to it.
Note that the standard presentation of Young tableaux as two-dimensional arrays of integers has
a significant disadvantage with respect to the computational cost and memory usage. That is because in this case the Sch\"{u}tzenberger transformation requires renumbering of all boxes in a tableau.

Let us consider the implemented algorithm for the case of 2D Young tableaux.
During the execution of the algorithm, the coordinates of boxes of a tableau are processed consequently.
At the beginning, the first box of a tableau with coordinates (0,0) is removed from a tableau and assigned as an \textit{active} box. Next boxes are being processed one by one. An active box is being inserted to a tableau at the moment when a neighbour top or neighbour right box was processed. As a next step, a neighbour box is assigned as an active box and being removed from a tableau and so on. At the same time, all non-neighbour boxes remain in place. The algorithm stops when all the boxes in a tableau are processed. 

The pseudocode of the two-dimensional Sch\"{u}tzenberger transformation algorithm is provided in Listing 1. The three-dimensional variant of the algorithm is very similar to it and is not presented here. Note that $actX, actY$ are the coordinates of the current active box, $tab$ is a Young tableau.

\smallskip

\begin{lstlisting}[language=C++, caption={Sch\"{u}tzenberger transformation on 2D Young tableaux\\~}]
schutz(tab)
{
  // remove the first box from tableau
  tab.remove(0,0)
  // initialize the coordinates of an active box
  actX = 0; actY = 0; 

  // for each (x,y) from tableau
  for (it = tab.begin(); it != tab.end(); ++it) 
  {
    // if the box (x,y) from a tableau is one
    // of two possible neighbours to an active box
    if ((x == actX + 1) && (y == actY)) || 
       ((x == actX)     && (y == actY + 1))
    {
      // remove it from tableau
      tab.remove(x,y);
      // add an active box to a tableau
      tab.insert(actX,actY);    
      // update an active box
      actX = x; actY = y;
    }
  }
  // return the coordinates of the removed box
  return (actX, actY);     
}
\end{lstlisting}

Note that during the execution of the algorithm, most of the boxes of a tableau remain in place. The number of modifications in a tableau is equal to the length of its Sch\"{u}tzenberger path. 
Another advantage of this approach is that after necessary modifications it can be easily implemented on a Young graph of any dimension and on any other graded graphs.

The implementation of the shape-preserving Sch\"{u}tzenberger transformation is almost the same. The only difference is that at the end of the algorithm the coordinates of the removed box should be added to a tableau. The reverse shape-preserving Sch\"{u}tzenberger transformation can also be implemented in a similar way.

The algorithm of the randomized Sch\"{u}tzenberger transformation in 2D case, written in pseudocode, is shown in Listing 2.

\smallskip

\begin{lstlisting}[language=C++, caption={Randomized Sch\"{u}tzenberger transformation on 2D Young tableaux\\~}]
schutz_rnd(tab, addlast)
{
  // if the beginning of tableau includes
  // Young diagram {2,1}
  if ( (tab[1] == (0,1)) && (tab[2] == (1,0)) || 
             (tab[1] == (1,0)) && (tab[2] == (0,1)) )
  { 
    // the path is chosen randomly from one  
    // of two options:
    // (0,0) => (0,1) => (1,0) or
    // (0,0) => (1,0) => (0,1)
    if (rand(0|1) == 0)
    {
      tab[1] = (0,1); tab[2] = (1,0);
    }
    else
    {
      tab[1] = (1,0); tab[2] = (0,1);
    }
  }
  // the classical Schutzenberger transformation function  
  (actX,actY) = schutz(tab);

  // add the last box of a nerve to a tableau
  if (addlast)
    tab.add(actX,actY);   
}
\end{lstlisting}

Unlike the shape-preserving Sch\"{u}tzenberger transformation, in the randomized modification the beginning part of the initial tableau is being either modified or left unchanged with a probability of 1/2.
Below, we will consider a three-dimensional variant of the randomized Sch\"{u}tzenberger transformation which differs from 2D case in some insignificant details.

\section{The applications of the randomized Sch\"{u}tzenberger's jeu de taquin}

\subsection{Calculation of the co-transition probabilities of the central process on 3D Young graph} \label{subs:cotrans}
The calculation of co-transition probabilities is a significant computational problem which appears in the study of Markov processes on graded graphs.
Let us consider a Markov process on the Young graph. \textit{The co-transition probability} $\widetilde{p}(\lambda_{n-1} \nearrow \lambda_{n})$ is a conditional probability that the path to diagram $\lambda_n$ passed through diagram $\lambda_{n-1}$, i.e. the ratio of dimensions $\frac{\dim{\lambda_{n-1}}}{\dim{\lambda_n}}$. There exists a hook-length formula which allows to calculate the exact dimensions of 2D Young diagrams. Then, the co-transition probability of a central process can be calculated by dividing the dimension of  $\lambda_{n-1}$ by dimension of  $\lambda_{n}$. Although, for the three-dimensional case such a formula is unknown.

The approach to estimate co-transition probabilities is to generate uniformly-distributed random paths to the given diagram and calculate the frequencies of passes through different incoming branches. Taking into account the law of large numbers, with increasing number of generations these frequencies will tend to unknown co-transition probabilities.

Unfortunately, the convergence rate of frequency ratios to co-transition probabilities is quite slow. For example, consider a pair of 3D Young diagrams of size 59 
$$
\lambda_{59}\!=\!\!\{\{7, 5, 4, 3, 2, 1, 1\}, \{5, 4, 3, 2, 1\}, \{4, 3, 2, 1\}, \{3, 2, 1\}, \{2, 1\}, \{1\}, \{1\}\}
$$
 and of size 60 
$$
\lambda_{60}\!=\!\!\{\{7, 5, 4, 3, 2, 2, 1\}, \{5, 4, 3, 2, 1\}, \{4, 3, 2, 1\}, \{3, 2, 1\}, \{2, 1\},
\{1\}, \{1\}\}.
$$
 The estimation of co-transition probability between them $E_{\lambda_{59} \nearrow \lambda_{60}} = 0.079$ which differs in the fourth decimal place from the ratio of exact dimensions of $\lambda_{59}$ and $\lambda_{60}$: $C_{\lambda_{59} \nearrow \lambda_{60}} = 0.079498$. However, we do not know other ways to calculate these probabilities even approximately.

\subsection{Calculation of dimensions of three-dimensional Young diagrams} 
Due to lack of analogue of the above-mentioned hook-length formula in a 3D case, it is not easy to calculate the dimensions of 3D Young diagrams.
The dimension of an arbitrary Young diagram can be calculated recurrently by dividing the dimension of a previous diagram by the corresponding co-transition probability: 
$$
\dim(\lambda_n) = \frac{\dim(\lambda_{n-1})}{\widetilde{p}(\lambda_{n-1} \nearrow \lambda_n)}.
$$
 So, these dimensions  can be estimated using the technique described in Sec.~\ref{subs:cotrans}. The exact   dimensions can be calculated by exhaustive search but such a calculation is possible for 3D diagrams of sizes only up to several tens of boxes.

Note that it is also possible to implement an algorithm for calculating the ratio of dimensions of diagrams nested inside one another. To do this, it is necessary to compare the products of co-transition probabilities of paths going from the intersection of these diagrams to each diagram. 
\newpage
\subsection{Three-dimensional Young diagrams of large dimensions}

The 2D and 3D Young diagrams with large and maximum dimensions are of special interest in asymptotic combinatorics. 
The problem of detection of 2D Young diagrams with large dimensions was studied in many works \cite{verpavl, ius15, pdmi15, knots, pdmi16, mcs16}. In the three-dimensional case the estimations of dimensions of diagrams are currently unknown.
The maximum dimensions of 3D Young diagrams grow extremely fast and it is very hard to calculate their exact values. At the moment, this problem can be solved only by exhaustive search.
Unfortunately, this approach gives only the first several tens of diagrams with maximum dimensions.
Table 1 contains the maximum dimensions of 3D Young diagrams of sizes up to 33 and the corresponding plane partitions.

One of possible approaches to find the diagrams with large dimensions is based on building a so-called greedy sequence of diagrams \cite{ius15, knots, pdmi15, mcs16}.
A greedy sequence of 2D diagrams is a deterministic sequence where the next diagram is obtained from the previous one by adding the box with the maximum Plancherel probability. In 3D case we can add the box with lowest possible co-transition probability considering that for 2D case these principles are equivalent.
Generally, a greedy sequence can be started from an arbitrary diagram. However, here we consider the sequence started from the root of the Young graph.

We have calculated exact dimensions of the first 65 diagrams of the greedy sequence. But we do not know whether the dimensions of diagrams from this sequence of sizes larger than 33 are maximum ones. We expect that many of them have maximum dimensions. The information about the first 65 terms of the greedy sequence is listed in Table 2. It includes the estimations of co-transition probabilities done using randomized Sch\"{u}tzenberger transformation, exact co-transition probabilities calculated by dividing exact dimensions of neighbour diagrams, the ratio of estimated and exact co-transition probabilities, and the ratio of estimated and exact dimensions.

\begin{table}[h]
\begin{center}
\captionof{table}{The maximum dimensions of plane partitions of size up to 33.}
{\footnotesize\begin{tabular}{|l|l|l|}
\hline
\textbf{Size} & \textbf{Dimension} & \textbf{Plane partition} \\
\hline
1	 & 	1	 & 	\{1\}\\
2	 & 	1	 & 	\{1\}, \{1\}\\
3	 & 	2	 & 	\{1 1\}, \{1\}\\
4	 & 	6	 & 	\{2 1\}, \{1\}\\
5	 & 	12	 & 	\{2 1\}, \{1\}, \{1\}\\
6	 & 	30	 & 	\{2 1 1\}, \{1\}, \{1\}\\
7	 & 	96	 & 	\{2 1 1\}, \{1 1\}, \{1\}\\
8	 & 	336	 & 	\{3 1 1\}, \{1 1\}, \{1\}\\
9	 & 	1540	 & 	\{3 1 1\}, \{2 1\}, \{1\}\\
10	 & 	8640	 & 	\{3 2 1\}, \{2 1\}, \{1\}\\
11	 & 	33372	 & 	\{3 2 1\}, \{2 1\}, \{1\}, \{1\}\\
12	 & 	142380	 & 	\{3 2 1 1\}, \{2 1\}, \{1\}, \{1\}\\
13	 & 	665280	 & 	\{4 2 1 1\}, \{2 1\}, \{1\}, \{1\}\\
14	 & 	2849536	 & 	\{3 2 1 1\}, \{2 1 1\}, \{1 1\}, \{1\}\\
15	 & 	15639552	 & 	\{4 2 1 1\}, \{2 1 1\}, \{1 1\}, \{1\}\\
16	 & 	80923008	 & 	\{4 2 1 1\}, \{2 1 1\}, \{2 1\}, \{1\}\\
17	 & 	544659648	 & 	\{4 2 1 1\}, \{3 1 1\}, \{2 1\}, \{1\}\\
18	 & 	3299672408	 & 	\{4 2 2 1\}, \{3 1 1\}, \{2 1\}, \{1\}\\
19	 & 	27402967200	 & 	\{4 3 2 1\}, \{3 1 1\}, \{2 1\}, \{1\}\\
20	 & 	230747045760	 & 	\{4 3 2 1\}, \{3 2 1\}, \{2 1\}, \{1\}\\
21	 & 	1553327915040	 & 	\{4 3 2 1\}, \{3 2 1\}, \{2 1\}, \{1\}, \{1\}\\
22	 & 	11012504995800	 & 	\{4 3 2 1 1\}, \{3 2 1\}, \{2 1\}, \{1\}, \{1\}\\
23	 & 	82028814137424	 & 	\{5 3 2 1 1\}, \{3 2 1\}, \{2 1\}, \{1\}, \{1\}\\
24	 & 	491203179370484	 & 	\{5 3 2 1 1\}, \{3 2 1\}, \{2 1\}, \{1 1\}, \{1\}\\
25	 & 	3290489409458592	 & 	\{5 3 2 1 1\}, \{3 2 1\}, \{2 1\}, \{2 1\}, \{1\}\\
26       &      26095216322563200        &      \{5 3 2 1 1\}, \{3 2 1 1\}, \{2 1 1\}, \{1 1\}, \{1\}\\
27       &      194868626458488668       &      \{5 3 2 1 1\}, \{3 2 1 1\}, \{2 1 1\}, \{2 1\}, \{1\}\\
28       &      1524692991397340664      &      \{5 3 2 1 1\}, \{3 2 1 1\}, \{2 1 1\}, \{2 1\}, \{1\}, \{1\}\\
29       &      13746015864155603608     &      \{5 3 2 1 1\}, \{4 2 1 1\}, \{3 1 1\}, \{2 1\}, \{1\}\\
30       &      118306078695096215552    &      \{5 3 2 1 1\}, \{4 2 1 1\}, \{3 1 1\}, \{2 1\}, \{1\}, \{1\}\\
31       &      1061302053614351456088   &      \{5 4 2 2 1\}, \{4 2 1 1\}, \{3 1 1\}, \{2 1\}, \{1\}\\
32       &      11607738064362975821328  &      \{5 4 3 2 1\}, \{4 2 1 1\}, \{3 1 1\}, \{2 1\}, \{1\}\\
33       &      111121303575872975022096 &      \{5 4 3 2 1\}, \{4 2 1 1\}, \{3 1 1\}, \{2 1\}, \{1\}, \{1\}\\
\hline
\end{tabular}}
\end{center}
\end{table}

\begin{table}
\begin{center}
\captionof{table}{The first 65 terms of the greedy sequence of 3D Young diagrams.}
{\footnotesize\begin{tabular}{|l|l|l|l|l|}
\hline
\textbf{Size} & \textbf{Estimation of} & \textbf{Exact co-transition} & \textbf{Ratio of co-transition} & \textbf{Ratio of exact} 
\\
 & \textbf{co-transition } & \textbf{probability} & \textbf{probabilities} & \textbf{dimensions}\\
& \textbf{probability} &&& \\
\hline
1 & 1.000000 & 1.000000 & 1.000000 & 1.000000\\
2 & 1.000000 & 1.000000 & 1.000000 & 1.000000\\
3 & 0.494000 & 0.500000 & 0.988000 & 0.988000\\
4 & 0.336400 & 0.333333 & 1.009200 & 0.997090\\
5 & 0.501300 & 0.500000 & 1.002600 & 0.999682\\
6 & 0.398900 & 0.400000 & 0.997250 & 0.996933\\
7 & 0.311000 & 0.312500 & 0.995200 & 0.992148\\
8 & 0.286700 & 0.285714 & 1.003450 & 0.995571\\
9 & 0.218100 & 0.218182 & 0.999625 & 0.995197\\
10 & 0.177300 & 0.178241 & 0.994722 & 0.989945\\
11 & 0.258500 & 0.258900 & 0.998456 & 0.988416\\
12 & 0.234000 & 0.234387 & 0.998350 & 0.986785\\
13 & 0.214000 & 0.214015 & 0.999929 & 0.986715\\
14 & 0.244400 & 0.243717 & 1.002804 & 0.989482\\
15 & 0.171500 & 0.174540 & 0.982583 & 0.972248\\
16 & 0.193500 & 0.193265 & 1.001218 & 0.973432\\
17 & 0.147200 & 0.148575 & 0.990743 & 0.964421\\
18 & 0.165300 & 0.165065 & 1.001425 & 0.965795\\
19 & 0.118700 & 0.120413 & 0.985774 & 0.952056\\
20 & 0.117900 & 0.118758 & 0.992779 & 0.945181\\
21 & 0.148500 & 0.148550 & 0.999663 & 0.944862\\
22 & 0.139800 & 0.141051 & 0.991129 & 0.936480\\
23 & 0.134000 & 0.134252 & 0.998125 & 0.934725\\
24 & 0.166900 & 0.166996 & 0.999427 & 0.934189\\
25 & 0.149900 & 0.149280 & 1.004155 & 0.938071\\
26 & 0.142200 & 0.141968 & 1.001636 & 0.939605\\
27 & 0.140300 & 0.140507 & 0.998527 & 0.938221\\
28 & 0.105200 & 0.108191 & 0.972356 & 0.912285\\
29 & 0.121100 & 0.123072 & 0.983975 & 0.897666\\
30 & 0.118500 & 0.118282 & 1.001846 & 0.899323\\
31 & 0.119200 & 0.119676 & 0.996022 & 0.895745\\
32 & 0.095700 & 0.092418 & 1.035517 & 0.927560\\
33 & 0.102400 & 0.102332 & 1.000666 & 0.928178\\
34 & 0.084100 & 0.084177 & 0.999086 & 0.927329\\
35 & 0.098300 & 0.097902 & 1.004063 & 0.931097\\
36 & 0.105300 & 0.103923 & 1.013254 & 0.943438\\
37 & 0.080800 & 0.079967 & 1.010416 & 0.953265\\
38 & 0.071400 & 0.070766 & 1.008953 & 0.961799\\
39 & 0.112400 & 0.113609 & 0.989361 & 0.951566\\
40 & 0.102000 & 0.103360 & 0.986841 & 0.939044\\
41 & 0.096900 & 0.096553 & 1.003593 & 0.942418\\
42 & 0.103600 & 0.103421 & 1.001734 & 0.944052\\
43 & 0.093600 & 0.095077 & 0.984470 & 0.929392\\
44 & 0.089100 & 0.089197 & 0.998908 & 0.928377\\
45 & 0.096400 & 0.095600 & 1.008368 & 0.936145\\
46 & 0.087900 & 0.087233 & 1.007646 & 0.943303\\
47 & 0.084400 & 0.082854 & 1.018658 & 0.960903\\
48 & 0.088500 & 0.090526 & 0.977617 & 0.939396\\
49 & 0.068400 & 0.068140 & 1.003811 & 0.942976\\
50 & 0.087400 & 0.086093 & 1.015178 & 0.957289\\
51 & 0.065800 & 0.064634 & 1.018034 & 0.974552\\
52 & 0.082300 & 0.081673 & 1.007679 & 0.982036\\
53 & 0.061200 & 0.061647 & 0.992753 & 0.974919\\
54 & 0.082300 & 0.082246 & 1.000659 & 0.975562\\
55 & 0.077700 & 0.078665 & 0.987731 & 0.963593\\
56 & 0.067500 & 0.068341 & 0.987690 & 0.951731\\
57 & 0.072900 & 0.073244 & 0.995305 & 0.947262\\
58 & 0.058600 & 0.059935 & 0.977719 & 0.926157\\
59 & 0.053300 & 0.053921 & 0.988483 & 0.915490\\
60 & 0.079000 & 0.079498 & 0.993739 & 0.909758\\
61 & 0.073700 & 0.073063 & 1.008718 & 0.917690\\
62 & 0.066200 & 0.067502 & 0.980715 & 0.899992\\
63 & 0.075800 & 0.075000 & 1.010664 & 0.909590\\
64 & 0.067300 & 0.069331 & 0.970711 & 0.882949\\
65 & 0.064400 & 0.063948 & 1.007074 & 0.889195\\
\hline
\end{tabular}}
\end{center}
\end{table}

\pagebreak
Due to the large values of exact dimensions, we did not include them in the Table 2. 
For example, the dimension of the diagram of size 65 
\{\{7, 5, 4, 3, 2, 2, 1, 1\}, \{5, 4, 3, 2, 1, 1\}, \{4, 3, 2, 1\}, \{3, 2, 1\}, \{2, 1\}, \{2, 1\}, \{1\},\{1\}\}
 is equal to 
$$11784492700515017182137999923695941374020209092205536828352.$$

There are 25 diagrams with maximum possible dimensions listed in Table 1 among the first 33 diagrams of this sequence. That is because there are no greedy sequences consisting only of diagrams with maximum dimensions. But it is possible to find the diagrams with larger dimensions than non-maximum ones from greedy sequences by applying some strategies discussed in \cite{ius15, pdmi15, knots}.

\section{Pseudo-Plancherel process}
The methods of generation of diagrams with large dimensions based on construction of greedy sequences and co-transition probabilities produce diagrams of sizes up to several thousands of boxes. However, some asymptotic problems require sequences of diagrams consisting of millions of  boxes. We use a special Markov process called pseudo-Plancherel process which was introduced in~\cite{terent} to generate diagrams of such sizes. 
Although this process is not central, computer experiments show that the probability ratios of paths between a pair of large Young diagrams are close to one. In this sense, the pseudo-Plancherel process can be considered as an asymptotically central one.

The formula for transition probabilities of this process is a direct generalization of a transition probabilities formula for 2D Plancherel process \cite{kerov93}. The transition probabilities of pseudo-Plancherel process are defined by the weight function which depends upon the lengths of 3D hooks in a diagram $\lambda$:

\begin{equation}
\label{eq:3dhook}
w(\lambda,x,y,z)\!=\! \!  \prod_{i=0}^{x-1}\! \frac{h(\lambda,i,y,z)}{h(\lambda,i,y,z)\!+\!1}\! \prod_{j=0}^{y-1}\! \frac{h(\lambda,x,j,z)}{h(\lambda,x,j,z)\!+\!1} \!\prod_{k=0}^{z-1} \!\frac{h(\lambda,x,y,k)}{h(\lambda,x,y,k)\!+\!1},
\end{equation}
where $x,y,z$ are coordinates of the added box, $h(\lambda,x,y,z)$ is the length of 3D hook of the $\lambda$ with the vertex at $x,y,z$. A 3D hook of a box ($x,y,z$) is a set of boxes $x', y', z'$ which satisfy the following condition: 
$$
\left[
	\begin{gathered} 
		(x' \geq x) \wedge (y' = y) \wedge (z' = z) \\
		(x' = x) \wedge (y' \geq y) \wedge (z' = z) \\
		(x' = x) \wedge (y' = y) \wedge (z' \geq z)
	\end{gathered}  
\right.
$$
The length of a 3D hook is the number of boxes in a hook.
The Pseudo-Plancherel transition probabilities are calculated by dividing the weights obtained in~\eqref{eq:3dhook} by the sum of all these values.

\section{Numerical experiments} \label{sec:numexp}

\subsection{Investigation of the randomized Sch\"{u}tzenberger transformation}
In this section we study the distribution of 2D and 3D Young tableaux produced by the random generator based on the randomized Sch\"{u}tzenberger transformation. The goal of these experiments is to examine how close this distribution is to the uniform distribution.

\subsubsection{The distribution of ends of nerves}
The first numerical experiment is devoted to the Sch\"{u}tzenberger transformation on 3D Young graph. 

\begin{figure}[ht]
	\centering
	\includegraphics[scale=1.45]{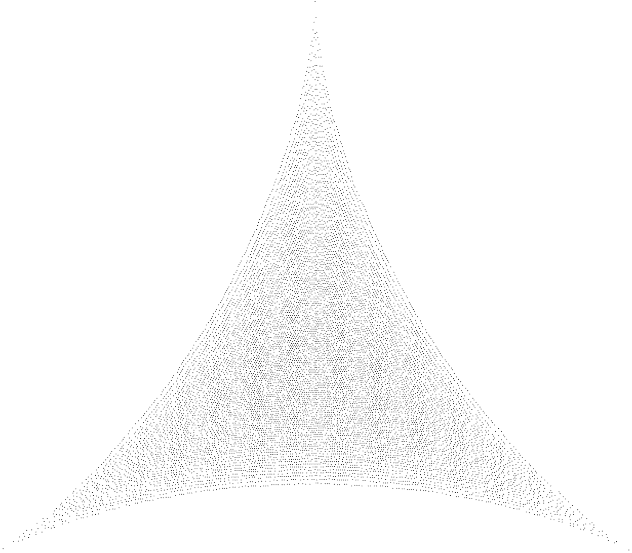}
	\caption{The distribution of coordinates of last boxes of Sch\"{u}tzenberger paths of 3D Young tableaux.}
	\label{fig:3D_lasty}
	\end{figure}

A random pseudo-Plancherel three-dimensional Young tableau of size $10^6$ was generated. Then, the randomized Sch\"{u}tzenberger transformation was applied iteratively to this tableau. For each iteration we compute the coordinates of the last boxes of Sch\"{u}tzenberger paths on the front of the corresponding Young tableau. The distribution of the coordinates obtained in this experiment is shown in Fig.~\ref{fig:3D_lasty}.

We emphasize that just after 66577 iterations, the Sch\"{u}tzenberger path endings in each of 11392 possible positions were generated at least once. This indicates that even for relatively small number of randomized Sch\"{u}tzenberger iterations, one of the necessary uniformity conditions for the generated tableaux is satisfied.

\subsubsection{The frequency histograms and statistical testing of the random sequences of Young tableaux}

\begin{figure}[!h]
        \centering
        \subfigure[]{
        \includegraphics[scale=0.345]{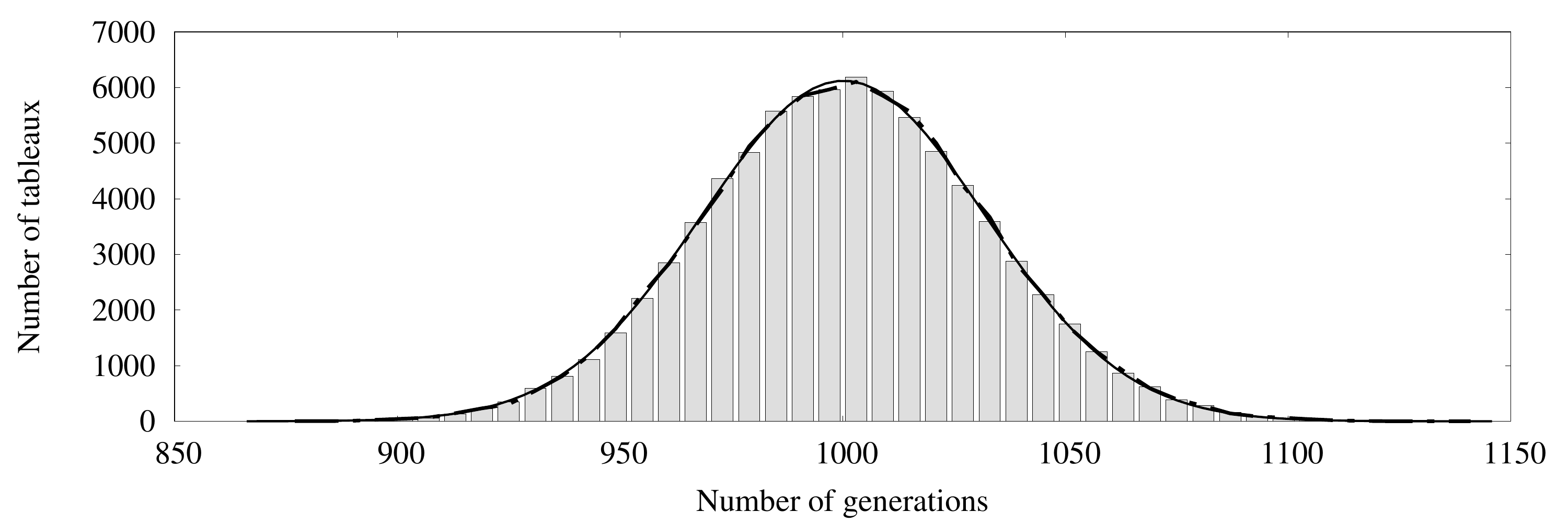}
        }
        \subfigure[]{
        \includegraphics[scale=0.345]{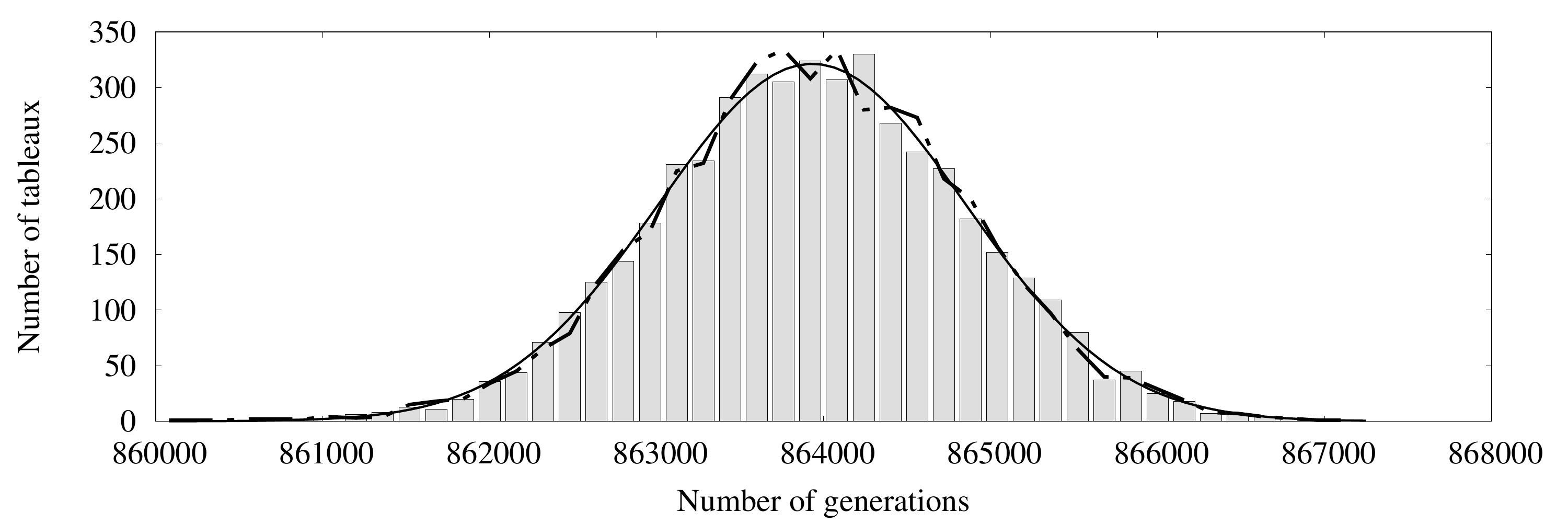}
        }
    \caption{Frequency histograms of random Young tableaux of the shape (a) $\lambda_1$ (2D), (b) $\lambda_2$ (3D). The dashed lines are the frequency histograms of uniformly distributed tableaux. The solid lines are Gaussian fits.}
        \label{fig:randomiz_hists}
\end{figure}

Let us consider the two-dimensional Young diagram $\lambda_1$ of size 15 which corresponds to the partition $\{4,4,3,3,1\}$. The dimension of this diagram is 81081. In our experiments we have iteratively applied the randomized Sch\"{u}tzenberger transformation to Young tableaux of that shape. 81081000 iterations were performed in total. The frequencies of generations of each tableau were recorded. Fig.~\ref{fig:randomiz_hists} (a) demonstrates how many tableaux were generated a certain number of times. 

As we can see, the mean of this graph is close to 1000 which is equal to the number of iterations divided by the dimension of the diagram. It agrees with the conjecture of the uniformity of this Young tableaux distribution. This histogram fits a normal distribution of standard deviation $\sigma= 31.7$ well which also indicates that the resulting Young tableaux are distributed almost uniformly.

Similarly, Fig.~\ref{fig:randomiz_hists} (b) presents the distribution of 3D Young tableaux of the shape $\lambda_2$ $\{\{2,2,1\},\{2,1\},\{1\},\{1\}\}$. $4\cdot10^9$ iterations were performed in total and the resulting histogram also fits a normal distribution well with mean $\mu = 863932$ and standard deviation $\sigma = 918.126$. The solid lines in Fig.~\ref{fig:randomiz_hists} represent Gaussian fits.

 An additional experiment was conducted in order to test the statistical hypothesis about the uniform distribution of a sequence of Young tableaux generated by iterations of the randomized 
Sch\"{u}tzenberger transformation. We used the standard C++ pseudo-random generator to produce 81081000 integers from the range [1..81081] and $4\cdot10^9$ integers from the range [1..4630]. The numbers were interpreted as random uniformly distributed numbers of tableaux. Then, the frequency histograms of these random integers were constructed in the same way as described above. The resulting histograms are depicted as dashed lines in Fig.~\ref{fig:randomiz_hists}. It can be easily seen that these histograms fit histograms of tableaux produced by randomized Sch\"{u}tzenberger transformation very well.
 This indicates a very good agreement between the test results and the hypothesis of a uniform distribution of the generated tableaux. 

\subsection{Calculation of normalized dimensions} \label{sec:calc_dim}

In view of the exponential growth of dimensions of 3D diagrams, it is convenient to study the asymptotic properties of dimensions using a certain normalization. In this work we used the following formula for normalized dimension of a three-dimensional Young diagram $\lambda$ of size $n$:
\begin{equation}
\label{eq:3dnormdim}
c(\lambda) = \frac{-\ln{\dim{\lambda}} + \frac{2}{3} \ln{n!}}{n^{\frac{2}{3}}},
\end{equation}
where $\dim{\lambda}$ is the dimension of $\lambda$.

Such a normalization was chosen because the leading term in the asymptotic expansion $\ln{\dim{\lambda}}=\frac{2}{3} \ln{n!}$.
Note that the formula~\eqref{eq:3dnormdim} is a generalization of normalized dimension formula of 2D Young diagram introduced in  \cite{verker85}.

In order to find the 3D Young diagrams with large dimensions, we consider a greedy sequence of Young diagrams  
described in Sec.4.3. Fig. 5 shows  how the normalized dimensions change in the greedy sequence.

	\begin{figure}[ht]
	\centering	\includegraphics[scale=0.55]{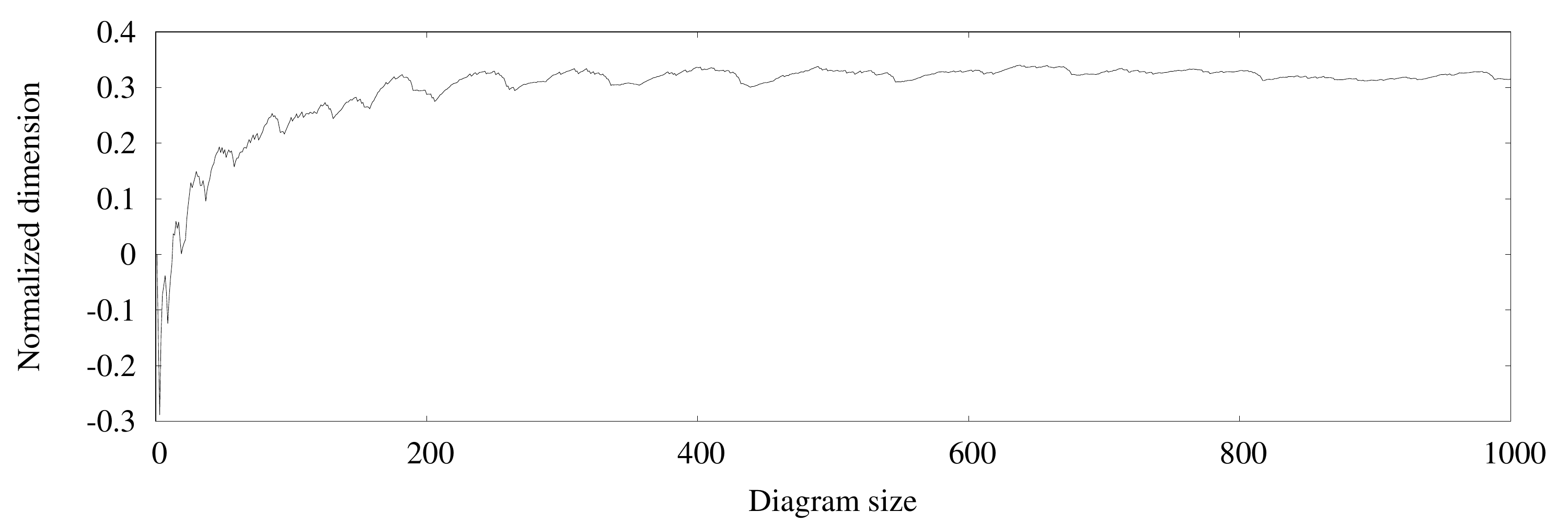}
	\caption{Approximated normalized dimensions of 3D Young diagrams of the greedy sequence.}
	\label{fig:3D_mincoprob}
	\end{figure}

The constructed sequence contains diagrams of sizes up to 1000. Alternatively, the longer greedy sequences can be efficiently produced by pseudo-Plancherel process instead. Computer experiments show that the dimensions of such diagrams are smaller.

\section{Conclusion}

The generator of random Young tableaux of a fixed shape based on the randomization of Sch\"{u}tzenberger's jeu de taquin was introduced. The analysis of results of our numerical experiments show that the distribution of tableaux produced by this generator is close to uniform. We propose a method of calculating the co-transition probabilities of the unknown central process on the three-dimensional Young graph using this generator. The algorithm for constructing sequences of 3D Young diagrams with large dimensions based on calculation of co-transition probabilities was presented. It was used to produce a sequence of  1000  3D Young diagrams with large dimensions. Normalized dimensions of diagrams from this sequence were calculated. It is expected that the asymptotic properties of this sequence is similar to the sequence of 3D Young diagrams with maximum dimensions.
\vspace{-.7em}
\subsection*{Acknowledgements}
The authors are grateful to Anatoly Vershik for numerous discussions on the issues raised in the article and to Vadim Smolensky for his comments that greatly improved the manuscript.

\end{document}